\documentclass[11pt]{amsart}
\usepackage{amsmath,amsfonts,amssymb,eucal}
\usepackage[all]{xy}

\begin{document}

\newcommand{\ie}{\textit{i}.\textit{e}.\,}
\newcommand{\eg}{\textit{e}.\textit{g}.\,}
\newcommand{\cf}{{\textit{cf}.\,}}

\newcommand{\ab}{{\rm ab}}
\newcommand{\alg}{{\rm alg}}
\newcommand{\Aut}{{\rm Aut}}
\newcommand{\bk}{{\bar{k}}}
\newcommand{\cC}{{\mathcal{C}}}
\newcommand{\C}{{\mathbb{C}}}
\newcommand{\End}{{\rm End}}
\newcommand{\sE}{{\sf E}}
\newcommand{\F}{{\mathbb{F}}}
\newcommand{\FG}{{\sf FG}}
\newcommand{\G}{{\mathbb{G}}}
\newcommand{\Gal}{{\rm Gal}}
\newcommand{\hp}{{\hat{p}}}
\newcommand{\Hom}{{\rm Hom}}
\newcommand{\iso}{{\sf iso}}
\newcommand{\la}{{\langle}}
\newcommand{\ra}{{\rangle}}
\newcommand{\Ksz}{{\rm Ksz}}
\newcommand{\Laz}{{\mathbb{L}}}
\newcommand{\LT}{{\sf LT}}
\newcommand{\m}{{\mathfrak{m}}}
\newcommand{\mor}{{\rm mor}}
\newcommand{\oh}{{\mathcal{O}}}
\newcommand{\oD}{{\mathfrak{o}_D}}
\newcommand{\oL}{{\mathfrak{o}_L}}
\newcommand{\ord}{{\rm ord}}
\newcommand{\cO}{{\mathcal{O}}}
\newcommand{\Qp}{{\mathbb{Q}_p}}
\newcommand{\Qq}{{\mathbb{Q}_q}}
\newcommand{\Q}{{\mathbb{Q}}}
\newcommand{\QZ}{{\mathbb{Q}/\mathbb{Z}}}
\newcommand{\rH}{{\rm H}}
\newcommand{\Spec}{{\rm Spec}}
\newcommand{\Spf}{{\rm Spf}}
\newcommand{\sF}{{\sf F}}
\newcommand{\sk}{{\sf k}}
\newcommand{\sK}{{\sf K}}
\newcommand{\THH}{{\rm THH}}
\newcommand{\zK}{{{\sf K}{\rm sz}}}
\newcommand{\trab}{{\rm trab}}
\newcommand{\tv}{{\tilde{v}}}
\newcommand{\U}{{\rm U}}
\newcommand{\Z}{{\mathbb{Z}}}

\title{Operations on integral lifts of $K(n)$}

\author{Jack Morava}

\address{Department of Mathematics, The Johns Hopkins University,
Baltimore, Maryland 21218}

\email{jack@math.jhu.edu}

\date{21 June 2019}

\begin{abstract}{This very rough sketch is a sequel to [27]; it presents evidence that operations 
on lifts of the functors $K(n)$ to cohomology theories with values in modules over valuation rings 
$\oL$ of local number fields, indexed by Lubin-Tate groups of such fields, are extensions of the 
groups of automorphisms of the associated group laws, by the exterior algebras on the normal bundle to 
the orbits of the group laws in the space of lifts.}\end{abstract}

\maketitle \bigskip

{\bf Introduction}

{\bf 0.1} In a symmetric monoidal category, \eg of schemes or structured spectra,
the morphisms defining an action of a monoid $M$ on an object $X$ can be presented 
as a cosimplicial object; for example [24] if $M = M\U$ is the Thom spectrum for complex
cobordism (\ie the universal complex-oriented $S^0$-algebra), then 
\[
\xymatrix{
( S^0 \ar@{.>}[r] & ) \; M\U \ar@<1ex>[r]\ar@<-1ex>[r] &  M\U \wedge_{S^0} M\U
\ar@<1ex>[r]\ar[r]\ar@<-1ex>[r] & \cdots}
\]
is a kind of $M\U$-free Adams-Mahowald-Novikov resolution of $S^0$. Its homotopy groups 
define a cosimplicial commutative algebra resolution
\[
\xymatrix{
\pi_*S^0 \ar@{.>}[r] & \pi_* M\U = M\U_* \ar@<1ex>[r]\ar@<-1ex>[r] & \pi_* ( M\U \wedge
_{S^0} M\U) = M\U_*M\U \ar@<1ex>[r]\ar[r]\ar@<-1ex>[r] & \cdots}
\]
of the stable homotopy algebra.\begin{footnote}{The terms in this display are
graded, but it is convenient to regard them as $\Z_2$-graded comodules over the 
multiplicative groupscheme $\G_m = \Z[t_0^{\pm 1}]$, with coaction $M_{2k} \ni x \mapsto 
x \otimes t_0^k$, thus providing an excuse for often supressing this grading.}\end{footnote}
Regarding these algebras as affine schemes over $\Spec \; \Z$, this diagram becomes a 
presentation for a groupoid-scheme 
\[
\xymatrix{
\Spec \; M\U_*M\U \ar@<1ex>[r]\ar@<-1ex>[r] & \Spec \; M\U_* }
\]
which, by work of Quillen [29], can be identified with a moduli stack for one-dimensional 
commutative formal groups. 

Using a great deal of work on Lubin-Tate spectra by others, we construct in \S 2.3.2 below, certain 
($p$-adically complete, where $p>3$) $A_\infty$ periodic $M\U$-algebra spectra $\sK(L)$, indexed 
by Lubin-Tate formal group laws $\LT_L$ for local number fields $L \supset \Qp$, Galois of degree 
$[L:\Qp] = n$ with valuation rings $\oL$. These spectra have homotopy groups 
\[
\pi_*\sK(L) = \sK(L)_* \cong \oL_*[v^{\pm 1}] 
\]
($|v| = 2$), and in \S 3 we present a conjectural description of an associated groupoid-scheme
\[
\xymatrix{
\Spec \; \sK(L)_*\sK(L) \ar@<1ex>[r]\ar@<-1ex>[r] & \Spec \; \sK(L)_* }
\]
of homological co-operations in terms of the isotropy or stabilizer groups of $\LT_L$, 
as objects in the Quillen-Lazard moduli stack. These automorphism groups are by now well-understood,
almost classical in local arithmetic geometry, and the first section below summarizes some
of that knowledge; it will serve as a model for our applications to algebraic topology.

Perhaps the point of this paper is to explain that, in spite of the notation, our construction
of the spectra $\sK(L)$ is {\bf not} functorial in $L$; this note is instead a plea for a natural
construction. The third section below contains preliminary results toward an identification of the 
endomorphisms or (co)operations of their associated cohomology theories, and argues that these
have close connections with the Weil group of $L$ [25]: or, more precisely, with the Galois 
group of a maximal totally ramified abelian extension $L^\trab$ of $L$, over $\Qp$. Our partial 
results can perhaps be read as evidence toward an interpretation of the spectra $\sK(L)$ as 
something like a $K$-theory spectrum associated to the (topological, perfectoid) completion 
$L^\infty$ of $L^\trab$ [28]. For example, our $\sK(\Qp)$ can be (non-canonically) identified 
with the $p$-adically completed algebraic $K$-theory spectrum of the completion $\Qp^\infty$ 
of the field of $p$-power roots of unity over $\Qp$, and thus with the $p$-adic completion of 
Atiyah's topological $K$-theory of $\C$.

To return to the organization of this paper: its second section uses the theory of highly 
structured spectra to define, following the original work of Sullivan and Baas [4,38,39], the spectra
$\sK(L)$ as Koszul quotients of spectra $\sE(\Phi_L)$ associated to Lubin-Tate formal 
group laws [14,31]. The resulting constructions are integral lifts of the `extraordinary' spectra
$K(n)$ [44], in that smashing with a mod $p$ Moore spectrum defines natural isomorphisms
\[
\sK(L)_*(X \wedge M(p,0)) \cong K(n)_*(X,\F_p) \otimes \oL/p\oL \;,
\]
where $\oL/p\oL \cong \F_q[\pi]/(\pi^e)$ (with $n = ef$ and $q = p^f$, see \S 1.5). For example, 
if $L$ is unramified then $e=1,\; q = p^n$, and the mod $p$ reduction of $\sK(L)$ agrees with 
$K(n) \otimes \F_q$. In some sense the $K(n)$ are indexed by the finite fields, while the $\sK(L)$ 
are indexed by finite Galois extensions of $\Qp$. \bigskip

{\bf \S I Notation and recollections}

{\bf 1.1} If $A$ is a commutative ring, let $\FG(A) \subset A[[X,Y]]$ be the set of power series 
$F(X,Y) = X + Y + \dots$ satisfying the standard axioms for a commutative formal group law over 
$A$, and let $\Gamma(A) \subset A[[T]]$ be the group of invertible power series $t(T) = t_0T + 
\dots$ (\ie with $t_0 \in A^\times$) under composition; then the group $\Gamma$ acts on the set 
$\FG$ by 
\[
\Gamma(A) \times \FG(A) \ni t,F \mapsto F^t(X,Y) = t^{-1}(F(t(X),t(Y))) \in \FG(A) \;.
\]
Both $\Gamma$ and $\FG$ are co-representable functors: $\FG(A) \cong \Hom_\alg(\Laz,A)$, where 
Lazard's ring $\Laz$ is polynomial over $\Z$, and $\Gamma(A) \cong \Hom_\alg(S,A)$, where $S = 
t_0^{-1}\Z[t_i]_{i \geq 0}$ is a Hopf algebra with coproduct
\[
(\Delta t)(T) = (t \otimes 1)((1 \otimes t)(T)) \in (S \otimes_\Z S)[[T]] \;.
\]
Yoneda's lemma then implies the existence of a coproduct homomorphism 
\[
\psi : \Laz \to \Laz \otimes_\Z S 
\]
of rings, corepresenting the group action. These rings are implicitly graded by the coaction
of the multiplicative subgroup $\G_m \subset \Spec \; S$.  

{\bf 1.2} A group action $\alpha : G \times X \to X$ in (Sets) defines a groupoid 
\[
\xymatrix{
[X//G] :  G \times X \ar@<1ex>[r]^-s \ar@<-1ex>[r]_-t & X}
\]
with $X$ as set of objects, $G \times X$ as set of morphisms, and $s(g,x) = x, \; t(g,x) = 
\alpha(g,x)$ as source and target maps. The usual convention in algebraic topology 
regards $\Laz \otimes_\Z S$ as a two-sided $\Laz$-algebra, with the obvious structure on the left,
and a right $\Laz$-algebra structure 
\[
\xymatrix{
(\Laz \otimes_\Z S) \otimes_\Z \Laz \ar[r]^-{1 \otimes \psi}  & (\Laz \otimes_\Z S) \otimes_\Laz 
(\Laz \otimes_\Z S) \ar[r] & (\Laz \otimes_\Z S) \;;}
\]
this is what's meant by saying that 
\[
\xymatrix{
\Laz \ar@<1ex>[r]^-{\eta_L} \ar@<-1ex>[r]_-{\eta_R} & \Laz \otimes_\Z S }
\]
is a Hopf algebroid.

Following Grothendieck and Segal, a category $\cC$ with set $\cC[0]$ of objects and $\cC[1]$ 
of morphisms can be presented as a simplicial set 
\[
\xymatrix{
\cC[0] & \ar@<1ex>[l] \ar@<-1ex>[l] \cC[1] & \ar@<1ex>[l] \ar[l] \ar@<-1ex>[l] \cC[1] 
\times_{\cC[0]} \cC[1] & \cdots }
\]
(where $X \times_Z Y$ denotes the fiber product or equalizer of two maps $X,Y \to Z$). In the case
of a group action as above, this is isomorphic to a simplicial object 
\[
\xymatrix{
X  & \ar@<1ex>[l] \ar@<-1ex>[l] G \times X  & \ar@<1ex>[l] \ar[l] \ar@<-1ex>[l] G \times G \times X
& \cdots }
\]
which can alternatively be regarded as a bar construction. The functor $A \mapsto 
[\FG(A)//\Gamma(A)]$ thus defines a simplicial scheme: the moduli stack of one-dimensional formal 
groups. 

{\bf 1.3} A homomorphism $A \to B$ of commutative rings defines an extension of scalars 
map 
\[
\FG(A) \ni F \mapsto F \otimes_A B \in \FG(B) \;. 
\]
{\bf Definition} {\it $[\iso(F)](B)$ is the groupoid with the orbit 
\[
\cO_\Gamma(B)(F) = \{ (F \otimes_A B)^g \:|\: g \in \Gamma(B) \}
\]
(of $F \otimes_A B$ under coordinate changes) as its set of objects, and
\[
\mor_{\iso_F(B)}(G,G') = \{ h \in \Gamma(B) \:|\: G^h = G' \}
\]
as (iso)morphisms of $G$ with $G'$; thus} 
\[
[\iso(F)](B) = [\cO_\Gamma(B)(F)//\Gamma(B)] \;.
\]
This groupoid maps fully and faithfully to its skeleton (which has one object) and the group 
$\Aut_B(F) \subset \Gamma(B)$ (of automorphisms of $F \otimes_A B$ as a formal group
law over $B$) as its morphisms. The homomorphism $F : \Laz \to A$ classifying $F$ thus defines
a Hopf $A$ - algebroid
\[
\xymatrix{
[\iso(F)] : A \ar@<1ex>[r] \ar@<-1ex>[r] & A \otimes_\Laz (\Laz \otimes_\Z S) \otimes_\Laz A }
\]
equivalent to a simplicial groupoid-scheme $(A - \alg) \ni B \mapsto [\iso(F)](B)$ over $\Spec 
\; A$. 

{\bf 1.4} It is nonstandard, but it will be convenient below to write $q = p^n$ and let $\Qq$ denote 
the quotient field of the ring $W(\F_q)$ of Witt vectors, \ie the degree $n$ unramified extension of 
$\Qp$. Following Ravenel [30 \S 5.1.13], a Lubin-Tate group law $\LT_\Qq$ for this field can be defined 
over the $p$-adic integers $\Z_p$ by Honda's logarithm 
\[
\log_\Qq(T) = \sum_{k \geq 0} p^{-k} T^{p^{nk}} \;;
\]
this has, as its mod $p$ reduction, a formal group law $\rH(n)$ of height $n$ over $\F_p$, 
associated to the cohomology theory $K(n)$. The resulting left and right $\F_p$-algebra structures on 
\[
\F_p \otimes_{\rH(n)} (\Laz \otimes_\Z S) \otimes_{\rH(n)} \F_p = 
C(\oD^\times,\F_q)^{\Gal(\F_q/\F_p)-{\rm inv}} = \Sigma(n) 
\]
coincide, representing $[\iso(\rH(n))]$ by the algebra of functions on a certain pro-algebraic 
group scheme over $\Spec \; \F_p$. 

In more detail [24], a finite field $k = \F_q$ has a local domain $W(k)$ of Witt vectors, with maximal ideal 
generated by $p$ and a canonical isomorphism $W(k)/pW(k) \to k$; its quotient field $W(k) \otimes_\Z \Q = 
\Qq$ is the extension of $\Qp$ obtained by lifting the roots $\F_q^\times$ of unity to $\Qp$. This 
construction is functorial, and a generator $\sigma$ of the cyclic group $\Gal(\Qq/\Qp) \cong 
\Gal(\F_q/\F_p)$ sends a root $\omega$ of unity to $\sigma(\omega) = \omega^p$. Let 
\[
D = \Qq \la \sF \ra/(\sF^n = p)
\]
be the noncommutative division algebra obtained from $\Qq$ by adjoining an $n^{th}$ root $\sF$ of 
$p$ satisfying, for any $a \in \Qq$, the relation $\sigma(a) \cdot \sF = \sF \cdot a$. The 
valuation on $\Qq$ (normalized so $\ord(p) = 1$) extends to $D$ to define a semidirect product extension
\[
\xymatrix{
1 \ar[r] & \oD^\times \ar[r] & D^\times \ar[r]^-\ord & \frac{1}{n} \Z \ar[r] & 0 }
\]
with a generator of the infinite cyclic group on the right acting on an element $u$ of the compact 
kernel $\oD^\times$ as $\sF$-conjugation. This kernel thus acquires an action of the cyclic group of 
order $n$, which may be identified with $\Gal(\F_q/\F_p)$, making $\oD^\times$ the group of points of 
a pro-\'etale groupscheme over $\F_p$. It is represented by the $\F_p$-algebra of (Galois equivariant, 
continous) $\F_q$-valued functions $h$ on $\oD^\times$ satisfying $\sigma(h(u)) = h(\sF u \sF^{-1})$. 
More concisely,
\[
[\iso(\rH(n))] \simeq [*//\oD^\times]
\]
as groupoid-valued functors.

{\bf 1.5} Similar results [21] hold for Lubin-Tate groups of local number fields (\ie extensions 
$L$ of $\Qp$ with $[L:\Qp] = n < \infty$); I will assume here that this extension is Galois. 
Such a field has a local valuation ring $\oL$ with finite residue field $k_L \cong \F_q$, where 
now $q = p^f$; moreover $L$ contains a maximal unramified extension $L_0 = W(k_L) \otimes_\Z \Q
\supset \Qp$, such that $[L:L_0] = e = f^{-1}n$. The maximal ideal $\m_L =(\pi) \subset \oL$ is 
principal, and we will choose a generator $\pi$; it  satisfies some Eisenstein equation 
\[
E_L(\pi) = \pi^e + \sum_{0 \leq i < e} e_i \pi^i = 0
\]
with $\ord(e_i) > 0$ and $\ord(e_0) = 1$. Lubin and Tate construct from this data, a 
formal group law $\LT_L$ over $L$ (with logarithm $\log_L$ and exponential $\exp_L$) such that 
\[
\oL \ni a \mapsto [a]_L(T) = \exp_L(a \cdot \log_L(T)) \in \End{\oL}(\LT_L)
\]
is a ring isomorphism. The reduction $\Phi_L$ over the residue field of $\LT_L$ is independent 
of the choices.

We will sometimes write $X +_L Y = \LT_L(X,Y)$. There is some degree of choice in the construction
of the lift $\LT_L$ of $\Phi_L$: 
\[
\log_L(T) = \sum_{i \geq 0} \pi^{-k} T^{q^k}
\]
is one possibility, and $[\pi]_L(T) = \pi T +_L T^q$ defines another\begin{footnote}{The author 
believes Serre's account [36] of the Lubin-Tate construction to be effectively optimal.}\end{footnote}
but all constructions are isomorphic. Reduction modulo $\pi$ defines a monomorphism
\[
\oL \cong \End_\oL(\LT_L) \to \End_{\bar{k}_L}(\Phi_L)
\]
of rings, which embeds the units $\oL^\times$ as a maximal commutative subgroup (a torus of some
sort) in $\oD^\times$.
 
The associated simplicial scheme $[\iso(\LT_L)]$ over $\Spec \; \oL$ can then be defined, as above, 
by the Hopf $\oL$-algebroid
\[
\xymatrix{
\oL \ar@<1ex>[r] \ar@<-1ex>[r] & \oL \otimes_{\LT_L} (\Laz \otimes_\Z S) \otimes_{\LT_L} \oL \;;}
\]
over the generic point of $\Spec \; \oL$ its fiber is the groupoid $[*//\oL^\times]$, while over
the closed point it is $[*//\oD^\times]$. Note that any degree $n$ extension of $\Qp$ embeds
in $D$ as a maximal commutative subfield, so the maximal toruses of $D^\times$ in some sense
parametrize Lubin-Tate groups of degree $n$ extensions of $\Qp$. The Weyl groups of these toruses 
are then Galois groups $\Gal(L/\Qp)$, and the normalizers of these toruses are essentially the Weil 
groups $W(L^\trab/\Qp)$ associated to maximal totally ramified abelian extensions of $L$ [41,43]; 
the (cohomology classes of the) group extensions defining them are the `fundamental classes' of local 
classfield theory.  

When $L = \Qq$ is unramified [27] we can assume that
\[
[p]_L(T) = pT +_L v^{q-1}T^q
\]
\ie that (a graded version of) $\LT_\Qq$ is $p$-typical, defined by a homomorphism 
\[
{\rm BP}_* = \Z_p[v_i]_{i \geq 1} \to \Z_p[v^{\pm 1}]
\]
sending Araki's [2] generators $v_i$ to 0 when $i \neq n$, and $v_n$ to $v^{q-1}$. \bigskip

{\bf \S II}

{\bf 2.0} To construct the spectra $\sK(L)$ we work\begin{footnote}{The typeface is intended to 
distinguish these constructions from Quillen's $K$-theory}\end{footnote} at a prime away from 6, in a 
symmetric monoidal category of $p$-adically complete spectra, \eg $S^0_\hp$ - modules. We will be 
concerned below with $K(n)$-local spectra, and we will $K(n)$-localize their smash products
[17]. Recall [33] that the Gaussian integer spectrum $S^0[(-1)^{1/2}]$ is not $E_\infty$:
the behavior of the stable homotopy category under arithmetic ramification seems potentially very
interesting.

Lurie's \'etale topology on the category of spectra [23 Def 7.5.1.4] defines commutative 
ring-spectra $S^0_{W(\F_q)}$ \'etale over $S^0_\hp$ (roughly, $W(\F_q)\otimes_{\Z_p} S^0_\hp)$. 
Schwede's Moore spectra (functorial away from 6 [34 \S II Rem. 6.44]) can then be used to define, 
for a valuation ring $\oL$ (free of rank $e$ over $W(k_L)$) a $p$-adic $A_\infty$ ring spectrum 
$S^0_\oL$ (roughly, $M(\oL,0) \otimes_{W(k_L)} S^0_{W(\F_q)}$) with 
\[
\pi_* S^0_\oL \cong \pi_* S^0 \otimes_\Z \oL \;.
\]
Following \S 1.5, the commutative W(k)-algebra
\[
\oL = \oplus_{0 \leq i \leq e-1} W(k) \cdot \pi^i
\]
(where $k = k_L$ for simplicity) is defined by classical structure constants
$m^{i,j}_l \in W(k), \; 0 \leq i,j,l \leq e-1$, such that
\[
\pi^i \cdot \pi^j = \sum_{0 \leq l \leq e-1} m^{i,j}_l \pi^l \;.
\]
Let $S^0_\oL$ denote the wedge sum $\bigvee_{0 \leq i \leq e-1} S^0_{W(k)} \cdot t^i$
(with $t^i$ a book-keeping indeterminate), and let
\[
S^0_\oL \times S^0_\oL \to S^0_\oL \wedge_{S^0_{W(k)}} S^0_\oL \to S^0_\oL
\]
be the morphism of $S^0_{W(k)}$-module spectra defined component-wise, as the composition
\[
\xymatrix{
S^0_{W(k)} \times S^0_{W(k)} \cdot t^i \times t^j \ar[r] &
S^0_{W(k)} \wedge_{S^0_{W(k)}} S^0_{W(k)} \cdot t^l \ar[r]^-{\cdot m^{i,j}_l} & S^0_{W(k)} t^l \;,}
\]
(where the final map is multiplication by the structure constant). This is the product map for a 
weak $S^0_{W(k)}$-algebra structure on $S^0_\oL$, \ie a kind of $H_\infty$ structure making
\[
\pi_*S^0_\oL = \pi_0 \Hom_{S^0_{W(k)}}(S^*_{W(k)},S^0_\oL) \cong \pi_*(S^0_\hp) \otimes_{\Z_p}
\oL
\]
as algebras. In particular we have a morphism
\[
\oL \times S^0_\oL \to S^0_\oL
\]
of $S^0_{W(k)}$-algebras, representing $\oL$-multiplication on $\pi_*S^0_\oL$, used in \S 2.3.2
below.

{\bf Remark} If $G$ is the Galois group of a finite extension $L$ of $\Qp$, a theorem
of Noether implies that its valuation ring $\oL$ is projective over the group ring $\Z_pG$
iff $L$ is tamely ramified; but work of Swan [40] implies, more generally, that the class of
$\oL$ in the Grothendieck group $G_0(\Z_pG)$ (defined by splitting short exact sequences)
is the image of a (not necessarily unique) class $[P]$ in $K_0(\Z_pG)$, perhaps analogous to
Wall's finiteness obstruction for CW complexes. Such Swan elements suggest constructing
analogs of Moore spectra for $\oL$ as representing objects for functors such as
\[
X \mapsto \pi_*(X \wedge G_+) \otimes_{\Z_pG} P := \pi_*(X;P) \; \dots
\]

{\bf 2.1.1} Work several mathematical generations deep [10,14,31,32\dots] associates to a 
one-dimensional formal group law $\Phi$, of finite height $n$ over a perfect field $k$ of 
characteristic $p > 0$, an $E_\infty \; p$-adic complex oriented $S^0_{W(k)}$-algebra spectrum 
$\sE(\Phi)$ with homotopy algebra 
\[
\pi*\sE(\Phi) \cong E(\Phi)_* \cong W(k)[[u_1,\dots,u_{n-1}]][v^{\pm 1}]]
\]
of formal power series, representing Lubin and Tate's functor [22] which sends a complete noetherian
local ring $A$ with residue field $k$ to the set (modulo isomorphisms which reduce to the identity 
over $k$) of lifts of $\Phi$ to $A$. We will sometimes take $v = 1$ to suppress the grading, and to 
simplify notation we may write $\sE_F$ for $\sE(\Phi)$ for a chosen lift $F$ of $\Phi$ to a local ring 
(\eg $W(k)$ or $\oL$) with residue field $k$; we may even write $\sE_{L_0}$ for the $S^0_{W(k_L)} =
S^0_{\oL_0}$-module spectrum $\sE(\LT_{L_0})$. Similarly $+_{E\Phi}$ or $+_{EF}$ may denote the 
associated formal group sum, and $\m_F = [[u_*]]$ may signify the `maximal ideal' of $E_{F*}$ over 
$W(k)$ or $\oL$. 

Lubin and Tate show that the (pro\'etale) group $\Aut_{\bk}(\Phi) \cong \oD^\times$ of automorphisms
of $\Phi \otimes_k \bk$, with its natural $\Gal(\bk/k)$-action, lifts to a (continuous but not 
smooth) action on $W(\bk) \otimes_{W(k)}E(\Phi)_*$; in particular, their theorem 3.1 shows that this 
action takes $W(\bk) \otimes_{W(k)} \m_{E(\Phi)}$ to itself. In the formalism of \S 1.3-4, this defines a 
formal groupoid-scheme of equivalence classes of lifts to Artin local rings, represented by a Hopf 
algebroid
\[
[\Spf \; E(\Phi)_*//\Aut_\bk(\Phi)] :\\
\xymatrix{  
E(\Phi)_* \ar@<1ex>[r]^-{\eta_L} \ar@<-1ex>[r]_-{\eta_R} & E(\Phi)_*\sE(\Phi) \cong H_{\Aut(\Phi)} \;;} 
\]
where $H_{\Aut(\Phi)}$ [17] is a Hopf algebra of Galois-equivariant continuous functions from $\Aut(\Phi)$ to
$E(\Phi)$. Note that the two (left and right) unit homomorphisms send $\m_E$ to $\m_E \hat{\otimes} 
H_{\Aut(\Phi)}$, and that this Hopf algebroid is equivalent to 
\[
\xymatrix{
E(\Phi)_* \ar@<1ex>[r]^-{\eta_L} \ar@<-1ex>[r]_-{\eta_R} & E(\Phi)_* \otimes_\Laz (\Laz \otimes_\Z) 
\otimes_\Laz E(\Phi)_* \;. }
\]

{\bf 2.1.2} If, for example, $\Phi_L/\F_p$ is the Lubin-Tate group law for $\LT_\Qq$ as in \S 1.3, 
we can take the $u_i$ to be Araki generators satisfying 
\[
[p]_E(T) = \sum_{E,i \geq o} v_i T^{p^i}
\]
(with $v_0 = p$); the classifying homomorphism from $M\U_*$ then sends $\C P_{q^k - 1}$ to 
\[
\prod_{1 \leq i \leq k} (1 - p^{q^i - 1})^{-1} \cdot (p^{-1}q^k)v^{q^k - 1} 
\]
and the remaining $\C P_l$ to 0 [17].

{\bf 2.2.1} A parallel (but even more venerable) line of research, leading to the modern theory of highly 
structured spectra, allows us to associate to the (by definition, regular) sequence $v_* = v_1,\dots,
v_{n-1}$ of elements of $E_{\Qq *}$, a choice
\[
v_i :S^{2(p^i - 1)} \to \sE_\Qq
\]
of representatives defining, by the construction of [12 V \S 1, \S 3.4], $p$-adic complex-oriented
$A_\infty$ ring-spectra
\[
\sK(\Qq) = \sE_\Qq/(v_*)
\]
with $\pi_*\sK(\Qq) \cong W(\F_q)[v^{\pm 1}]$, having $\LT_\Qq$ as formal group law.

{\bf Remark} When $n=1$ this construction recovers a model for Atiyah's $p$-adic completion [3] of 
complex topological $K$-theory, and when $n=2$ it defines a $p$-adic lift of Baker's supersingular 
elliptic cohomology [5]. Away from the prime 6, elliptic cohomology [13 \S 4.1 ex 4.2] has as 
coefficients, the polynomial ring of modular forms (generated [37 III \S 5.6.2] by the Eisenstein 
series $E_2,E_3$). A theorem of Deligne [19] identifies the modular form defined by the Eisenstein 
series $E_{p-1}$ and the (Hasse) parameter $v_1$, modulo $p$. This suggests a close relation between 
$p$-adic elliptic cohomology, mod $E_{p-1}$, with $\sK(\Q_{p^2})$; but understanding that would 
require an understanding of $E_{p-1}$ as a polynomial in $E_2,E_3$, which evidently depends on the 
prime $p$.

{\bf 2.2.2} More generally, a Lubin-Tate group law for a ramified local field of degree $n$ over 
$\Qp$ lifts its mod $\pi$ reduction to a homomorphism
\[
u_i \mapsto u^*_i : E(\Phi_L)_* = W(k_L)[[u_*]] \to \oL
\]
of local $W(k_L)$-algebras. 

{\bf Lemma} {\it Let $u^*_i \in \m_L = (\pi) \subset \oL, \; i \leq i \leq n-1$ be a sequence of elements
in the maximal ideal $\m_L$: then $u_i \mapsto u_i - u^*_i = v_i$ defines an isomorphism $\oL[[u_*]]
\cong \oL[[v_*]]$ of local $\oL$-algebras.}

[For if w$ \in \m_L$ and $\oL[[x]] \ni a(x) = \sum_{i \geq 1} a_i x^i$, then 
\[
a(x) = \sum_{k \geq 1} (\sum_{l \geq 1} \binom{k+l}{l} a_{k+l} w^l) \tilde{x}^k = \sum_{l \geq 1}
\tilde{a}_k \tilde{x}^k \;,
\]
where $\tilde{x} = x - w$. The argument for multiple variables is similar.] $\Box$

{\bf Definition} {\it The $A_\infty$ complex-oriented $S^0_\oL$-algebra spectrum 
\[
\sE_L = S^0_L \wedge_{S^0_{W(k_L)}} \sE(\Phi_{\LT_L})
\]
has 
\[
\pi_* \sE_L = E_{L*} = \oL[[v_*]][v^{\pm 1}]
\]
as algebra of homotopy groups, generated by a regular sequence of elements $v_i = u_i - u^*_i$
such that $v_i \mapsto 0$ specializes the modular lift to the chosen Lubin-Tate group law of $L$.}  

{\bf 2.3.1} The definition in \S 2.1 of $\sK(\Qq)$ uses the $E_\infty$ structure on 
$\sE(\Phi_{\LT_\Qq})$, which is not available for nontrivially ramified fields. This issue can be 
avoided by reorganizing the induction in [12] (which follows Sullivan and Baas, based on iterated 
cofibrations) as a computation of the spectral sequence for the homotopy groups of the geometric 
realization of a suitable simplicial (Koszul) spectrum [11 \S 17.5, ex 18.2; 35]:

An element $a$ of a commutative $k$-algebra $A$ defines\begin{footnote}{or, more generally, a homogeneous 
element of a graded commutative $A_*$; however this will be largely suppressed from our 
notation}\end{footnote} an (elementary) differential graded $A$-algebra 
\[
\Ksz_A(a) = (A[e]/(e^2), \; d_a (e)= a) \;.
\]
More generally, a $k$-module homomorphism $a_\star : k^m \to A$ defines the classical differential 
graded commutative algebra 
\[
\Ksz_A(a_\star) = \bigotimes_{A, 1 \leq i \leq m} \Ksz_A(a_i) \cong A \otimes_k \Lambda_k(e_i \:|\: 1 \leq i 
\leq m)
\]
(with an exterior algebra denoted by $\Lambda$, to reduce the multiplicity of things called $E$), 
and differential 
\[
d e_I = \sum_{1 \leq i \leq m} (-1)^{i+1} a_\star(e_i) \cdot e_{\hat{I}(i)} \;,
\]
where $e_I = \wedge_{l \in I} e_l$ is indexed by (totally ordered) subsets $I$ of $\{m\} = \{1,\dots,
m\}, \; \hat{I}(i)$ is obtained from $I$ by omitting its $i^{th}$ element, $e_I \wedge e_K$ is 0 if 
$I \cap K$ is nonempty and equals $\pm e_{I,K}$ if they are disjoint, with sign equal to that of 
the permutation putting $\{I,K\}$ in proper order.   

In the elementary case, if $a$ is not a 0-divisor in $A$, this defines an $A$-free resolution of the 
quotient algebra $A/(a)$, \ie of the cofiber of the map of $A$ to itself defined by 
$a$-multiplication. More generally, if $a_i$ is not a 0-divisor in the quotient ring 
$A/(a_1,\dots,a_{i-1})$ (\ie \:$a_*$ is a {\bf regular} sequence), this construction defines an 
$A$-free resolution of $A/(a_\star)$. 

In the context of commutative $S^0$-algebras or ring spectra, we can associate to morphisms
\[
a_i : S^{|a_i|} \to A 
\]
a semi-simplicial (\ie without degeneracy operators [42 \S 8.2.2]) $A$-algebra
\[
k \mapsto \zK_A(a_\star)[k] = \bigvee_{I \subset \{m\}, |I| = k} A \cdot e_I
\]
with $A$-module face operators $\partial_i e_I = \mu(a_i) \cdot e_{\hat{I}(i)}$, where
\[
\mu(a) : S^{|a|} A \to A 
\]
is defined by multiplication by $a$. The fat realization [11 \S 4.8, \S 18.2] 
$|\zK_A(a_\star)|$ of such a semisimplicial object is canonically filtered, leading to the
construction of a spectral sequence computing its homotopy groups, with $E_1$ page
\[
\pi_* \zK_A(a_\star) = \Ksz_{A_*}(a_\star) \Rightarrow |\zK_A(a_\star)|_* \;
\]
 
{\bf 2.3.2} When $a_\star$ is a regular sequence this complex is a resolution, and the spectral
sequence collapses to an isomorphism
\[
|\zK_A(a_\star)|_* \cong A_*/(a_\star) \;.
\]
Applying this as above yields a definition for $\sK(\Qq) = |\zK_{\Phi(\Qq)}(v_\star)|$ as an 
$S^0_{W(\F_q)}$-algebra spectrum, essentially equivalent to the construction of [12]. The ramified 
case is more delicate, because its building blocks are not $E_\infty$; we need a

{\bf Lemma} {\it The morphisms 
\[
\tv_i,\tv_j : S^*_\oL \sE_L \to \sE_L 
\]
($1 \leq i,j \leq n-1)$ defined by multiplication with 
\[
\tv_i = u^*_i \wedge_W 1_\sE - 1_L \wedge_W u_i : S^*_\oL \to S^0_\oL \wedge_{W(k_L)} \sE(\Phi_L) 
\; (\; = \sE_L)
\]
commute.}

{\bf Proof} Define a twist isomorphism
\[
\sE(\Phi_L) \wedge_{W(k_L)} \oL \to \oL \wedge_{W(k_L)} \sE(\Phi_L) 
\]
adjoint to the composition 
\[
\sE(\Phi_L) \to \Hom_{W(k_L)}(\oL,\oL) \wedge_{W(k_L)} \sE(\Phi_L) \to 
\Hom_{W(k_L)}(\oL,\oL \wedge_{W(k_L)} \sE(\Phi_L))
\]
of $S^0_{W(k_L)}$-module morphisms. Since both $S^0_L$ and $\sE(\Phi_L)$ are commutative 
$S^0_{W(k_L)}$-modules, we have (with some abbreviation)
\[
\tv_i \wedge_W \tv_j = (u^*_i \wedge_W 1_\sE - 1_L \wedge_W u_i) \wedge_W (u^*_j \wedge_W 1_\sE - 1_L 
\wedge_W u_i) = \cdots = \tv_j \wedge \tv_i \;.
\]

{\bf Proposition} 
\[
\sK(L) = |\zK_{\sE_L}(\tv_\star)|
\]
{\it is an $A_\infty \; S^0_\oL$-algebra spectrum with $\sK(L)_* \cong \oL[v^{\pm 1}]$, complex-oriented
by the morphism $M\U_* \to \sK(L)_*$ classifying the chosen Lubin-Tate group law for} $L$. $\Box$

{\bf Remark} The kernel of 
\[
\oL \otimes_\Z M\U_* \to E_{L*} \to \sK(L)_*
\]
is generated by an (infinite) regular sequence (which can be chosen to belong to $M\U_*$ in degree greater
than $2(p^n-1)$). The Koszul construction above, together with [12 V \S 3.2], leads to a construction for
a connective spectrum $\sk(L)$ with $\sK(L)$ as $v_n$-localization.\bigskip

{\bf \S III}

{\bf 3.1} The spectral sequence of a geometric realization, together with the Eilenberg-Moore/K\"unneth spectral
sequence for the smash product of module spectra provide some understanding of the bialgebra 
\[
(\sK(L) \wedge_{S^0_L} \sK(L)_* = \sK(L)_*\sK(L) \;.
\]
To begin, note that 
\[
(\sK(L) \wedge_{S^0_L} \sE_L)_* = |\zK_{\sE_L}(\tv_\star)|_*(\sE_L)
\]
is the $\sE(\Phi_L)_*$ homology of a filtered $\sE(\Phi_L)$-module spectrum, and that the $E_1$ page of the 
associated spectral sequence is the Koszul algebra 
\[
\Ksz_{(\sE_L \wedge_{S^0_L} \sE_L)_*}(\tv_\star) \;;
\]
but by [17], as in \S 2.1, $(\sE_L \wedge_{S^0_L} \sE_L)_* \cong H_{\Aut(\Phi_L)}$, with 
deformation parameters acting as left $\tv_*$-multiplication. This sequence is regular, so this is a 
resolution, and the spectral sequence collapses to an isomorphism 
\[
(\sK(L) \wedge_{S^0_L} \sE_L)_* \cong \sK(L)_* \otimes_{E_{L*}} H_{\Aut(\Phi_L)} \cong \oL 
\otimes_{W(k_L)} H_{\Aut(\Phi_L)} 
\]
of $\oL$-algebras. Now observe that 
\[
\sK(L) \wedge_{S^0_L} \sK(L) \simeq (\sK(L) \wedge_{S^0_L} \sE_L) \wedge_{\sE_L} \sK(L) 
\]
which is accessible via [12 IV Thm 6.4].

{\bf Proposition} {\it The K\"unneth spectral sequence collapses at $E_2$ to an isomorphism
\[
\sK(L)_*\sK(L) \cong \oL \otimes_{E_{L*}} H_{\Aut(\Phi_L)} \otimes_\oL \Lambda^*_\oL(\m_L/\m_L^2) \;,
\]
where the term on the right is the exterior algebra on the (free, of rank $n-1$) tangent $\oL$-module
to the space of deformations of $\LT_L$.} 

{\bf Proof} The $E_1$-page of this spectral sequence is again a Koszul algebra, now of the form
\[
\Ksz_{(\sK(L)_* \otimes_{\sE_{L*}} H_{\Aut(\Phi_L)})}(\eta(v_\star)) \;,
\]
where the images
\[
\eta(v_i) = \sum_\alpha v_{i,\alpha} \otimes g_{i,\alpha} \in \m_L H_{\Aut(\Phi_L)}
\]
of the generators $v_i$ under the right unit have coefficients $v_{i,\alpha}$ in the ideal $\m_L$ 
([22 Thm 3.1], see \S 2.1), and thus map to zero under $\sE_{L*} \to \sK(L)_*$. The homology of 
this DGA is therefore just its underlying graded algebra, which can be identified with the algebra
of Galois-equivariant functions from $\Aut(\Phi_L)$ to the exterior algebra on $\m_L/\m_L^2$. [Note 
that $\sE_{L*}\sE_L$ is analogous [by Kodaira-Spencer theory, \cf [15 ex 2.8.1, 20]] to an algebra 
of functions from $\Aut(\Phi_L)$ to the symmetric algebra on $\m_L/\m_L^2$.] $\Box$ 

{\bf 3.2} It seems reasonable to conjecture that this spectral sequence collapse implies an interpretation
of $\sK(L)_*\sK(L)$ as a Hopf algebroid of functions on a (super, ie nontrivially $\Z_2$-graded) groupoid
scheme, an extension of the automorphism group of $\LT_L$ by an exterior algebra of deformations parametrized
by its tangent space as a point in $\Spf \; \sE_{L*}$. However, the author feels that this and related 
questions (\eg the possible nontriviality of such extensions, the action of $\Gal(L/\Qp)$ and its previously 
mentioned relation to Weil groups, Massey product structures [1 \S 3.3], relations [16] with Azumaya 
algebras \dots) are best left to younger, more vigorous and reliable resarchers.

In particular: recent advances in the algebra of non-discretely valued fields suggest that the topological 
Hochschild homology of the perfectoid completion $L^\infty$ of the fields $L^\trab$ (\S 1.5) associates to
the $p$-adic completion of $B\QZ$ (regarded as an analog of $\C P^\infty$), a rigid analytic analog [28] of a
Lubin-Tate group for $L$. If the spectra $\sK(L)$ have a natural construction in terms of fields like 
$L^\infty$, one might hope for the existence of a generalized Chern character or cyclotomic-like trace, mapping
$\sk(L)$ Galois-equivariantly to $\THH(\oh_{L^\infty},\Z_p)$.

\bibliographystyle{amsplain}

\end{document}